\documentclass[12pt]{article}

\usepackage{amsmath,amsthm,amssymb,a0size,bm,indentfirst,txfonts}
\usepackage{stmaryrd,CJK,color,graphicx}
\usepackage[toc,page,title,titletoc,header]{appendix}
\usepackage{paralist}
\usepackage{epsfig}
\usepackage[colorlinks=true]{hyperref}
\hypersetup{citecolor=red}

\usepackage{remreset}

\allowdisplaybreaks

\newcommand{\bex}{\begin{eqnarray*}}
\newcommand{\eex}{\end{eqnarray*}}
\newcommand{\be}{\begin{eqnarray}}
\newcommand{\ee}{\end{eqnarray}}
\newcommand{\ba}{\begin{array}}
\newcommand{\ea}{\end{array}}
\newcommand{\bi}{\begin{itemize}}
\newcommand{\ei}{\end{itemize}}
\newcommand{\bn}{\begin{enumerate}}
\newcommand{\en}{\end{enumerate}}

\newcommand{\sex}[1]{\left(#1\right)}

\newcommand{\sez}[1]{\left[#1\right]}
\newcommand{\sed}[1]{\left\{#1\right\}}
\newcommand{\sev}[1]{\left|#1\right|}
\newcommand{\sen}[1]{\left\Vert#1\right\Vert}

\newcommand{\dps}{\displaystyle}
\newtheorem{thm}{\indent Theorem}
\newtheorem{cor}[thm]{\indent Corollary}
\newtheorem{lem}[thm]{\indent Lemma}

\newtheorem{defn}[thm]{\indent Definition}
\newtheorem{rem}[thm]{\indent Remark}

\newcommand{\bp}{\begin{proof}}
\newcommand{\ep}{\end{proof}}

\newcommand{\R}{\mathbb{R}}

\def\n{\nabla}
\def\p{\partial}
\def\lap{\Delta}

\def\bbu{\bm{u}}
\def\bbv{\bm{v}}

\def\bbC{\bm{C}}

\def\bbH{\bm{H}}
\def\bbL{\bm{L}}

\def\bbV{\bm{V}}

\def\bbcalV{\bm{\mathcal{V}}}

\def\bbphi{\bm{\phi}}

\makeatletter
\def\@eqnnum{{\normalfont \color{red} (\theequation)}}
\makeatother

\makeindex

\begin{document}
\title{A Serrin-type regularity criterion for the Navier-Stokes equations via one velocity component}

\author{Zujin Zhang
\thanks{Department of Mathematics,
Sun Yat-sen University, P.R. China 510275}
\thanks{Email: uia.china@gmail.com}}
\maketitle

\begin{abstract}
We study the Cauchy problem for the $3$D Navier-Stokes equations, and prove some scalaring-invariant regularity criteria involving only one velocity component.
\end{abstract}

\section{Introduction}
\label{sect:intro}
We consider the following incompressible Navier-Stokes equations in $\R^3$:
\be\label{eq:NSE}
\ba{cc}
\left.
\ba{cc}
\p_t\bbu-\nu\lap \bbu+\sex{\bbu\cdot\n}\bbu+\n p=0\\
\n\cdot\bbu=0
\ea
\right\}
&\mbox{in }(0,T)\times\R^3,\\
\bbu(0,x)=\bbu_0(x),
&\mbox{in }\R^3.
\ea
\ee
Here $\bbu=\sex{u_1,u_2,u_3}$ is the velocity, $p$ is a scalar pressure, $\nu>0$ is the kinematic viscosity, and $\bbu_0$ with $\n\cdot \bbu_0=0$ is the initial velocity.

The existence of a global weak solution
\bex
\bbu\in L^\infty(0,T;\bbL^2(\R^3))
\cap
L^2(0,T;\bbH^1(\R^3))
\eex
to \eqref{eq:NSE} has long been established by Leray \cite{Leray_34}, see also Hopf \cite{Hopf}. But the issue of regularity and uniqueness of $\bbu$ remains open. Initialed by Serrin \cite{Serrin_62,Serrin_63} and Prodi \cite{Prodi_59}, there has been a lot of literatures devoted to finding sufficient conditions to ensure $\bbu$ to be smooth, see, e.g., \cite{Beirao da Veiga_95},\cite{Beirao da Veiga_02}, \cite{Constantin_Fefferman_93}, \cite{Escauriaza_Seregin_Sverak_03}, \cite{Kim_10}, \cite{Neustupa_Penel}, \cite{Penel_Pokorny}, \cite{Zhou_05_direction vorticity}, \cite{Zhou_06_grad pressure}, \cite{Zhou_04_pressure}, \cite{Zhou_09_weighted reg}, \cite{Zhang_Xicheng_08}, \cite{Zhang_Chen_05}, and references therein.

Recently, many authors become interested in the regularity criteria involving only one velocity component, or its gradient, even though most of which are not scaling invariant. Let us track the progresses we made during the last decade.

For one component regularity, one is preferred to showing that the condition
\be\label{reg:u3}
u_3\in L^p(0,T;L^q(\R^3)),\ 2/p+3/q=\beta,\ 3/\beta<q\leq \infty,
\ee
with $\beta=1$ or
\be\label{reg:grad u3}
\n u_3\in L^r(0,T;L^s(\R^3)),\ 2/r+3/s=\gamma,\ 3/\gamma<q\leq\infty,
\ee
with $\gamma=2$, ensures that $\bbu\in \bbC^\infty((0,T]\times \R^3)$.

However, this is quite difficult to prove, and all results are with $\beta<1$ or $\gamma<2$,  to the authors' best knowledge. More precisely, $\beta$ (resp. $\gamma$) is first taken to be $1/2$ (resp. $3/2$) in \cite{Neustupa_Novotny_Penel_99} and \cite{Zhou_05_one component}. Then, using intricate decomposition of the pressure $p$, Kukavica and Ziane \cite{Kukavica_Ziane_06_one component} was able to show that \eqref{reg:u3} with $\beta=5/8$ or \eqref{reg:grad u3} with $\gamma=11/6$ is enough to ensure smoothness of $\bbu$. Later, Cao and Titi \cite{Cao_Titi_08_one component} extended $\alpha$ to be $3/2+2/3q$, by invoking muliplicative Sobolev imbedding inequality:
\be\label{ineq:mul Sobolev}
\sen{f}_6\leq C\sen{\n_h f}_2^{2/3}\sen{\p_3 f}_2^{1/3},
\ee
where $\n_h=\sex{\p_1,\p_2}$ is the horizontal gradient. Finally, Zhou and Pokorn\'y in \cite{Zhou_Pokorny_09_grad one component},\cite{Zhou_Pokorny_09_one component} give another contribution, which states that the condition
\eqref{reg:u3} with $\beta=3/4+1/2q$ or \eqref{reg:grad u3} with $\gamma=23/12$ entails the regularity of $\bbu$, although, there are some restrictions on $s$.

Interestingly enough, regularity in one direction is always scaling invariant, see \cite{Cao_2010_one component}, \cite{Kukavica_Ziane_07_one direction}.

The purpose of this paper is to make a further contribution
in this direction. Precisely, we have
\begin{thm}
\label{thm:main}
Let $\bbu_0\in \bbV$, and $\bbu$ be a weak solution to \eqref{eq:NSE} in $[0,T]$ with initial datum $\bbu_0$. If
\be\label{thm:main:reg}
u_3\in L^p(0,T;L^q(\R^3)),\ \p_3u_3\in L^r(0,T;L^s(\R^3)), \ee
with $1\leq p,q,r,s\leq \infty$, $0\leq \beta,\gamma<\infty$ satisfying
\be\label{thm:main:ass}
\left\{
\ba{ccc}
2/p+3/q=\beta,\ 2/r+3/s=\gamma\\
\dps{\sex{1-\frac{1}{s}}q=\frac{1/r+3/8}{3/8-1/p}=\frac{9/4-\gamma}{\beta-3/4}>1}\\
p<\infty\mbox{ or }r<\infty
\ea
\right.
\ee
then $\bbu$ is smooth in $[0,T]\times\R^3$.
\end{thm}
\begin{rem}
We make some comments on \eqref{thm:main:ass}. For the sake of corollaries followed, we write down $\eqref{thm:main:ass}_1$. In $\eqref{thm:main:ass}_2$, the first equality is due to some H\"older $/$ Young conjugates, see $\eqref{exponent}_{1,2}$ and the derivation of \eqref{final},
while the second one is some compatibility condition, see $\eqref{exponent}_3$. Finally, $\eqref{thm:main:ass}_3$ is to ensure the application of Gronwall inequality in \eqref{last 2nd}.
\end{rem}
If we choose $\beta=1$ or $\gamma=2$ in Theorem \ref{thm:main}, we have the following Serrin-type regularity criterionco.
\begin{cor}\label{cor:serrin}
Assume as in Theorem \ref{thm:main}. If
\bex
u_3\in L^p(0,T;L^q(\R^3)),\ \p_3u_3\in L^r(0,T;L^s(\R^3)),
\eex
with $1\leq p,r\leq\infty, 1\leq q,s\leq \infty$ satisfying either of the following conditions,
\bn
    \item
        \bex
        \left\{
        \ba{ccc}
        2/p+3/q=1,\ 2/r+3/s<2\\
        \dps{\sex{1-\frac{1}{s}}q=\frac{1/r+3/8}{3/8-1/p}
        =9-4\sex{\frac{2}{r}+\frac{3}{s}}}\\
        p<\infty\mbox{ or }r<\infty
        \ea
        \right.
        \eex
    \item
        \bex
        \left\{
        \ba{ccc}
        2/p+3/q<1,\ 2/r+3/s=2\\
        \dps{\sex{1-\frac{1}{s}}q=\frac{1/r+3/8}{3/8-1/p}
        =\frac{1}{4(2/p+3/q)-3}}\\
        p<\infty\mbox{ or }r<\infty
        \ea
        \right.
        \eex
\en
then $\bbu$ is smooth in $[0,T]\times \R^3$.
\end{cor}

Due to Definition \ref{defn:wf} in Sect. \ref{sect:pre}, $u_3\in L^\infty(0,T;L^2(\R^3))$, we may take $p=\infty$, $q=2$, $\gamma=3/4+3/2s$ in Theorem \ref{thm:main} to yield

\begin{cor}\cite{Cao_Titi_10_one entry reg} Assume as in Theorem \ref{thm:main}. Then the condition
\bex
\p_3u_3\in L^r(0,T;L^s(\R^3)),\ 2/r+3/s=3/4+3/2s,\ 1\leq r<\infty,
\eex
ensures that $\bbu\in \bbC^\infty([0,T]\times\R^3)$.
\end{cor}
The function space $\bbV$, the definition of a weak solution and other often-used notations will be given in Sect. \ref{sect:pre}.

We shall use method from \cite{Cao_Titi_10_one entry reg} and \cite{Zhou_Pokorny_09_one component}. See the details in Sect. \ref{sect:proof}.

\section{Prelimiaries}
\label{sect:pre}
We gather here some definitions, notations and
intricate (in)equalities.

The Lebesgue spaces $L^q(\R^3)$ is endowed with norm $\sen{\cdot}_q$, with its bold-face counterpart denotes the set of vector-valued functions, and we denote by $\sen{\cdot}_{p,q}$ the norm for anisotropic Lebesgue spaces $L^p(0,T;L^q(\R^3))$. The Sobolev spaces $W^{k,p}(\R^3)$ (resp. $H^k(\R^3)$) is equipped with the norm $\sen{\cdot}_{k,p}$ (resp. $\sen{\cdot}_{k,2}$).

Let $\bbC^\infty_c(\R^n)$ be the set of smooth vector-valued functions with compact support, we then define
\bex
\bbcalV=\sed{\bbv\in \bbC_c^\infty(\R^3);\ \n\cdot\bbv=0},
\eex
\bex
\bbH=\mbox{ the completion of }\bbcalV \mbox{ under the norm }\sen{\cdot}_2,
\eex
\bex
\bbV=\mbox{ the completion of }\bbcalV \mbox{ under the norm }\sen{\cdot}_{1,2}.
\eex

With these spaces at hand, we recall the weak formulation of \eqref{eq:NSE}, see \cite{Temam_77_NSE}.
\begin{defn}\label{defn:wf}
Let $\bbu_0\in \bbH$, $T>0$. A measurable vector-valued function $\bbu$ defined in $[0,T]\times\R^3$ is said to be a weak solution to \eqref{eq:NSE} if
\bn
    \item $\bbu\in L^\infty(0,T;\bbH)\cap L^2(0,T;\bbV)$;
    \item $\eqref{eq:NSE}_1$ holds in the sense of distributions,
        \be
        \label{defn:wf:eq:wf}
        & &\sex{\bbu(t),\bbphi(t)}
        +\nu\int_0^t\sex{\n\bbu(s),\n\bbphi(s)}ds\nonumber\\
        & &
        -\int_0^t\sex{\bbu,\p_t\bbphi+\sex{\bbu\cdot\n}\bbphi}ds
        =\sex{\bbu_0,\bbphi(0)},
        \ee
        for all $\bbphi\in \bbC_c^\infty([0,T)\times\R^3)$ with $\n\cdot\bbphi=0$.

        Here $\sex{\cdot,\cdot}$ is the scalar product in $\bbL^2(\R^3)$.
\en
\end{defn}

Since we are concerned with regularity criteria involving only one velocity component, the following lemma is quite important, see \cite{Kukavica_Ziane_07_one direction}.
\begin{lem}\label{lem:kukavica_ziane}
Assume that $\bbu\in \bbH^2(\R^3)$ is smooth and divergence free. Then
\be\label{lem:kukavica_ziane:eq:u_3}
& &\sum_{i,j=1}^2\int_{\R^3}u_i\p_iu_j\lap_hu_jdx\\
&=&\frac{1}{2}\sum_{i,j=1}^2
\int_{\R^3}\p_iu_j\p_iu_j\p_3u_3dx
-\int_{\R^3}\p_1u_1\p_2u_2\p_3u_3dx\nonumber\\
& &
+\int_{\R^3}\p_1u_2\p_2u_1\p_3u_3dx\nonumber\\
&=&-\sum_{i,j=1}^2\int_{\R^3}
\p_{3i}^2u_j\p_iu_ju_3dx
+\int_{\R^3}\p_{31}^2u_1\p_2u_2u_3dx\nonumber\\
& &
+\int_{\R^3}\p_1u_1\p_{32}^2u_2u_3dx\nonumber\\
& &-\int_{\R^3}\p_{31}^2u_2\p_2u_1u_3dx
-\int_{\R^3}\p_1u_2\p_{32}^2u_1u_3dx,\nonumber
\ee
where $\lap_h=\p_{11}^2+\p_{22}^2$ is the horizontal Laplacian.
\end{lem}

It is then immediate that by invoking the divergence free condition (see \cite{Cao_Titi_10_one entry reg})
\begin{lem}\label{lem:bdd dot horizontal lap}
\be\label{lem:bdd dot horizontal lap:ineq}
\sev{
\int_{\R^3}\sez{\sex{\bbu\cdot\n}\bbu}\cdot\lap_h \bbu dx}
\leq \int_{\R^3}
\sev{u_3}\cdot\sev{\n\bbu}\cdot \sev{\n\n_h\bbu}dx.
\ee
\end{lem}
We end this section by invoking some interpolation inequalities, see e.g., \cite{Zhou_05_one component}, and a simple revision of an inequality in \cite{Cao_Titi_10_one entry reg}.

\begin{lem}\label{lem:interpolation}
Let $f\in L^\infty(0,T;L^2(\R^3))\cap L^2(0,T;H^1(\R^3))$. Then
\bex
f\in L^a(0,T;L^b(\R^3)),\ 2/a+3/b\geq 3/2,\ 2\leq b\leq 6.
\eex
Moreover,
\be\label{lem:interpolation:ineq}
\sen{f}_{a,b}\leq C\sen{f}_{\infty,2}^{3/b-1/2}
\sen{\n f}_{2,2}^{3/2-3/b}.
\ee
\end{lem}

\begin{lem}\label{lem:trilinear estimates}
For $f$, $g$, $h$ $\in C_c^\infty(\R^3)$, we have
\be\label{lem:trilinear estimates:ineq}
& &\sev{\int_{\R^3}\ f\ g\ h\ dx_1dx_2dx_3}\nonumber\\
&\leq& C\sen{f}_q^\frac{\alpha-1}{\alpha}
\sen{\p_3f}_s^{1/\alpha}
\sen{g}_2^\frac{\alpha-2}{\alpha}
\sen{\p_1g}_2^{1/\alpha}
\sen{\p_2g}_2^{1/\alpha}
\sen{h}_2,
\ee
where
\bex
\alpha>2,\ 1\leq q,s\leq \infty,\ \frac{\alpha-1}{q}+\frac{1}{s}=1.
\eex
\end{lem}
\bp
\bex
& &\sev{\int_{\R^3}\ f\ g\ h\ dx_1dx_2dx_3}\\
&\leq&\int_{\R^2}
    \sez{
    \max_{x_3}\sev{f}
    \cdot
    \sex{\int_{\R}g^2dx_3}^{1/2}
    \cdot
    \sex{\int_{\R}h^2dx_3}^{1/2}
    }
    dx_1dx_2\\
& &\sex{\mbox{H\"older inequality}}\\
&\leq&\sez{
    \int_{\R^2}\sex{\max_{x_3}\sev{f}}^\alpha dx_1dx_2
    }
    \cdot
    \sez{
    \int_{\R^2}\sex{\int_{\R}g^2dx_3}
    ^\frac{\alpha}{\alpha-2}dx_1dx_2                 \cdot
    }^\frac{\alpha-2}{2\alpha}\\
& &\cdot
    \sex{\int_{\R^3}h^2dx_1dx_2dx_3}^{1/2}\ \sex{\mbox{H\"older inequality again}}\\
&\leq&C\sez{\int_{\R^3}\sev{f}^{\alpha-1}
    \sev{\p_3f}dx_1dx_2dx_3}^{1/\alpha}
    \cdot
    \sez{\int_{\R}
    \sex{\int_{\R^2}
    g^\frac{2\alpha}{\alpha-2}dx_1dx_2}
    ^\frac{\alpha-2}{\alpha}dx_3}^{1/2}
    \cdot
    \sen{h}_2\\
& &\sex{\mbox{Minkowski inequality}}\\
&\leq&C\sen{f}_q^\frac{\alpha-1}{\alpha}
\sen{\p_3 f}_s^{1/\alpha}
\sen{g}_2^\frac{\alpha-2}{\alpha}
\sen{\p_1 g}_2^{1/\alpha}
\sen{\p_2 g}_2^{1/\alpha}
\sen{h}_2\\
& &\sex{\mbox{H\"older, interpolation inequalities and }\eqref{ineq:mul Sobolev}}. \eex
\ep

\section{Proof of the main result}
\label{sect:proof}

In this section, we prove Theorem \ref{thm:main}.

{\bf Step I} Some reductions.

By the classical "weak = strong" type uniqueness theorem, we need only prove that
\be\label{goal}
\sen{\n \bbu}_{\infty,2}<\infty.
\ee
Due to \eqref{thm:main:ass}, we can take an $\alpha>2$ such that
\bex
\alpha-1=(1-1/s)q=\frac{1/r+3/8}{3/8-1/p}=\frac{9/4-\gamma}{\beta-3/4},
\eex
i.e.,
\be\label{exponent}
\left\{
\ba{ccc}
\dps{\frac{\alpha-1}{p}+\frac{1}{r}=\frac{3(\alpha-2)}{8}},\\
\dps{\frac{\alpha-1}{q}+\frac{1}{s}=1},\\
\dps{(\alpha-1)\beta+\gamma=\frac{3(\alpha-2)}{4}+3}.
\ea
\right.
\ee

{\bf Step II} $\sen{\n_h \bbu}_2$ estimates.

Taking the inner product of $\eqref{eq:NSE}_1$ with $-\lap_h \bbu$ in $\bbL^2(\R^3)$, we obtain
\be\label{grad hori u:process}
& &\frac{1}{2}\frac{d}{dt}\sen{\n_h\bbu}_2^2
+\nu\sen{\n_h\n \bbu}_2^2\nonumber\\
&=&\int_{\R^3}\sez{\sex{\bbu\cdot\n}\bbu}\cdot\lap_h\bbu dx\nonumber\\
&\leq&C\int_{\R^3}\sev{u_3}\cdot\sev{\n\bbu}\cdot\sev{\n\n_h\bbu}dx
    \ \sex{\mbox{by }\eqref{lem:bdd dot horizontal lap:ineq}}\nonumber\\
&\leq&C\sen{u_3}_q^\frac{\alpha-1}{\alpha}
    \sen{\p_3u_3}_s^{1/\alpha}
    \sen{\n \bbu}_2^\frac{\alpha-2}{\alpha}
    \sen{\n\n_h\bbu}_2^\frac{\alpha+2}{\alpha}
    \ \sex{\mbox{by }\eqref{lem:trilinear estimates:ineq}, \eqref{exponent}_2}\nonumber\\
&\leq&\frac{\nu}{2}\sen{\n\n_h\bbu}_2^2
    +C\sen{u_3}_q^\frac{2(\alpha-1)}{\alpha-2}
    \sen{\p_3u_3}_s^\frac{2}{\alpha-2}
    \sen{\n\bbu}_2^2.
\ee
Integrating \eqref{grad hori u:process}, we deduce
\be\label{grad hori u}
& &\sen{\n_h\bbu(t)}_2^2
+\nu\int_0^t\sen{\n_h\n\bbu(s)}_2^2ds\nonumber\\
&\leq&
    \sen{\n_h\bbu_0}_2^2+
    C\int_0^t\sen{u_3}_q^\frac{2(\alpha-1)}{\alpha-2}
    \sen{\p_3u_3}_s^\frac{2}{\alpha-2}
    \sen{\n\bbu(s)}_2^2ds,
\ee
for all $t\in [0,T]$.

{\bf Step III} $\sen{\n \bbu}_2$ estimates.

Taking inner product of $\eqref{eq:NSE}_1$ with $-\lap\bbu$ in $\bbL^2(\R^3)$, we gather, noticing \eqref{lem:bdd dot horizontal lap:ineq} and $\n\cdot\bbu=0$, that
\be\label{grad u:process}
\frac{1}{2}\frac{d}{dt}\sen{\n\bbu}_2^2
+\nu\sen{\lap\bbu}_2^2
&=&\int_{\R^3}\sez{\sex{\bbu\cdot\n}\bbu}\cdot\lap_h \bbu dx
    +\int_{\R^3}\sez{\sex{\bbu\cdot\n}\bbu}\p_{33}^2\bbu dx\nonumber\\
&\leq&C\int_{\R^3}\sez{\sev{u_3}\cdot\sev{\n\bbu}\cdot\sev{\n\n_h\bbu}
    +\sev{\n_h\bbu}\cdot\sev{\p_3\bbu}^2}dx\nonumber\\
&\equiv& J_1+J_2.
\ee
We estimate $J_1$ as in \eqref{grad hori u:process},
\bex
J_1\leq \frac{\nu}{2}\sen{\lap \bbu}_2^2
    +C\sen{u_3}_q^\frac{2(\alpha-1)}{\alpha-2}
    \sen{\p_3u_3}_s^\frac{2}{\alpha-2}
    \sen{\n \bbu}_2^2.
\eex
Meanwhile, using H\"older, interpolation inequalities and \eqref{ineq:mul Sobolev}, $J_2$ is dominated as
\bex
J_2&\leq&C\sen{\n_h\bbu}_2\sen{\n\bbu}_4^2\\
&\leq&C\sen{\n_h\bbu}_2\sen{\n\bbu}_2^{1/2}\sen{\n\bbu}_6^{3/2}\\
&\leq&C\sen{\n_h\bbu}_2\sen{\n\bbu}_2^{1/2}\sen{\n_h\n\bbu}_2\sen{\lap \bbu}_2^{1/2}.
\eex
Replacing these two last displaced inequalities into \eqref{grad u:process} yields
\be\label{grad u}
& &\frac{d}{dt}\sen{\n\bbu}_2^2
+\nu\sen{\lap\bbu}_2^2\nonumber\\
&\leq& C\sen{u_3}_q^\frac{2(\alpha-1)}{\alpha-2}
     \sen{\p_3u_3}_s^\frac{2}{\alpha-2}
     \sen{\n \bbu}_2^2
    +C\sen{\n_h\bbu}_2\sen{\n\bbu}_2^{1/2}\sen{\n_h\n\bbu}_2\sen{\lap \bbu}_2^{1/2}.
\ee
Integrating the above inequality, and invoking H\"older inequality, we have
\bex
& &\sen{\n\bbu(t)}_2^2
+\nu\int_0^t \sen{\lap \bbu}_2^2ds\\
&\leq&\sen{\n\bbu_0}_2^2
    +C\int_0^t \sen{u_3}_q^\frac{2(\alpha-1)}{\alpha-2} \sen{\p_3u_3}_s^\frac{2}{\alpha-2}\sen{\n\bbu}_2^2 d s
    +C\sup_{0\leq s\leq t}\sen{\n_h\bbu(s)}_2\\
& &\cdot\sex{\int_0^t\sen{\n\bbu(s)}_2^2ds}^{1/4}
    \cdot\sex{\int_0^t\sen{\n\n_h\bbu(s)}_2^2ds}^{1/2}
    \cdot\sex{\int_0^t\sen{\lap\bbu(s)}_2^2ds}^{1/4}.
\eex
Thanks to \eqref{grad hori u}, we get
\bex
& &\sen{\n\bbu(t)}_2^2
+\nu\int_0^t \sen{\lap \bbu}_2^2ds\\
&\leq&\sen{\n\bbu_0}_2^2
    +C\int_0^t \sen{u_3}_q^\frac{2(\alpha-1)}{\alpha-2}
    \sen{\p_3u_3}_s^\frac{2}{\alpha-2}\sen{\n\bbu}_2^2ds\\
& &
    +C\sez{
    \sen{\n_h\bbu_0}_2^2+
    C\int_0^t\sen{u_3}_q^\frac{2(\alpha-1)}{\alpha-2}
    \sen{\p_3u_3}_s^\frac{2}{\alpha-2}
    \sen{\n\bbu(s)}_2^2ds
    }
    \cdot\sex{\int_0^t\sen{\lap\bbu(s)}_2^2ds}^{1/4}.
\eex
By Young and H\"older inequalities, and the fact $\bbu\in L^2(0,T;V)$, we find
\be\label{last 2nd}
& &\sen{\n\bbu(t)}_2^2
+\frac{\nu}{2}\int_0^t\sen{\lap \bbu(s)}_2^2ds
    \nonumber\\
&\leq&C\sen{\n\bbu_0}_2^2
    +C\int_0^t\sen{u_3}_q^\frac{2(\alpha-1)}{\alpha-2}
    \sen{\p_3u_3}_s^\frac{2}{\alpha-2}
    \sen{\n\bbu(s)}_2^2ds\nonumber\\
& &
    +C\int_0^t\sen{u_3}_q^\frac{8(\alpha-1)}{3(\alpha-2)}
     \sen{\p_3u_3}_s^\frac{8}{3(\alpha-2)}
     \sen{\n\bbu(s)}_2^2ds.
\ee
Thanks to $\eqref{exponent}_1$ and Young inequality, we obtain
\be\label{final}
& &\sen{\n\bbu(t)}_2^2
+\frac{\nu}{2}\int_0^t\sen{\lap \bbu(s)}_2^2ds\nonumber\\
&\leq&C\sen{\n\bbu_0}_2^2
    +C\int_0^t
    \sex{\sen{u_3}_q^p
    +\sen{\p_3u_3}_s^r
    +1}
    \sen{\n\bbu(s)}_2^2ds.
\ee
Therefore, by Gronwall inequality and \eqref{thm:main:reg}, we have \eqref{goal} as desired. The proof is completed.


\begin{thebibliography}{00}

\bibitem{Beirao da Veiga_95}
H. Beir\~ao da Veiga,
\emph{A new regularity class for the Navier-Stokes equations in $\R^n$,}
Chinese Ann. Math. Ser. B
{\bf 16}
(1995),
407--412.

\bibitem{Beirao da Veiga_02}
H. Beir\~ao da Veiga, L.C. Berselli,
\emph{On the regularizing effect of the vorticity direction in incompressible viscous flows,}
Differential Integral Equations
{\bf 15}
(2002),
345--356.

\bibitem{Cao_2010_one component}
C.S. Cao,
\emph{Sufficient conditions for the regularity to the $3$D Navier-Stokes equations,}
Discrete Contin. Dyn. Syst.
{\bf 26}
(2010),
1141--1151.

\bibitem{Cao_Titi_10_one entry reg}
C.S. Cao, E.S. Titi,
\emph{Global regularity criterion for the $3$D Navier-Stokes equations involving one entry of the velocity gradient tensor,}
arXiv: 1005.4463 [math. AP] 25 May 2010.

\bibitem{Cao_Titi_08_one component}
C.S. Cao, E.S. Titi,
\emph{Regularity criteria for the three-dimensional Navier-Stokes equations,}
Indiana Univ. Math. J.
{\bf 57}
(2008),
2643--2661.

\bibitem{Constantin_Fefferman_93}
P. Constantin, C. Fefferman,
\emph{Direction of vorticity and the problem of global regularity for the Navier-Stokes equations,}
Indiana Univ. Math. J.
{\bf 42}
(1993),
775--789.

\bibitem{Escauriaza_Seregin_Sverak_03}
L. Escauriaza, G. Seregin, V. \v Sver\'ak,
\emph{Backward uniqueness for parabolic equations,}
Arch. Ration. Mech. Anal.
{\bf 169}
(2003),
147--157.

\bibitem{Fan_Jiang_Ni_08}
J.S. Fan, S. Jiang, G.X. Ni,
\emph{On regularity criteria for the $n$-dimensional Navier-Stokes equations in terms of the pressure,}
J. Differential Equations
{\bf 244}
(2008),
2963--2979.

\bibitem{Hopf}
E. Hopf,
\emph{\"Uer die Anfangswertaufgabe f\"ur die hydrodynamischen Grundgleichungen,}
Math. Nachr.
{\bf 4}
(1951),
213--231.

\bibitem{Kim_10}
J.M. Kim,
\emph{On regularity criteria of the Navier-Stokes equations in bounded domains,}
J. Math. Phys.
{\bf 51}
(2010),
053102.

\bibitem{Kukavica_Ziane_07_one direction}
I. Kukavica, M. Ziane,
\emph{Navier-Stokes equations with regularity in one direction.}
J. Math. Phys.
{\bf 48}
(2007), 065203.

\bibitem{Kukavica_Ziane_06_one component}
I. Kukavica, M. Ziane,
\emph{One component regularity for the Navier-Stokes equations.}
Nonlinearity
{\bf 19}
(2006),
453--469.

\bibitem{Leray_34}
J. Leray,
\emph{Sur le mouvement d'un liquide visqueux emplissant l'espace.}
Acta Math.
{\bf 63}
(1934),
193--248.

\bibitem{Neustupa_Novotny_Penel_99}
J. Neustupa, A. Novotn\'y, P. Penel,
\emph{An interior regularity of a weak solution to the Navier-Stokes equations in dependence on one component of velocity,}
Topics in Mathematical Fluid Mechanics,
Quaderni di Matematica, Dept. Math., Seconda University, Napoli, Caserta,
Vol. 10, pp. 163-183 (2002);
see also
\emph{A remark to interior regularity of a suitable weak solution to the Navier-Stokes equations,}
CIM Preprint No. 25, 1999.

\bibitem{Neustupa_Penel}
J. Neustupa, P. Penel,
\emph{Anisotropic and geometric criteria for interior regularity of weak solutions to the $3$D Navier¨CStokes equations,}
in Mathematical Fluid Mechanics (Recent Results and Open Problems), Advances in Mathematical Fluid Mechanics, edited by J. Neustupa, and P. Penel (Birkh\"auser, Basel-Boston-Berlin, (2001), pp. 239¨C267.

\bibitem{Penel_Pokorny}
P. Penel, M. Pokorn\'y,
\emph{Some new regularity criteria for the Navier¨CStokes equations containing the gradient of velocity,}
Appl. Math.
{\bf 49}
(2004),
483--493.

\bibitem{Prodi_59}
G. Prodi,
\emph{Un teorema di unicit\'a per le equazioni di Navier-Stokes,}
Ann. Mat. Pura Appl.
{\bf 48}
(1959),
173--182.

\bibitem{Serrin_62}
J. Serrin,
\emph{On the interior regularity of weak solutions of the Navier-Stokes equations,}
Arch. Rational Mech. Anal.
{\bf 9}
(1962),
187--191.

\bibitem{Serrin_63}
J. Serrin,
\emph{The initial value problems for the Navier-Stokes equations,}
in Nonlinear Problems, edited by R. E. Langer (University of Wisconsin Press, Madison, WI, (1963).

\bibitem{Temam_77_NSE}
R. Temam,
\emph{Navier-Stokes equations, Theory and numerical analysis,}
North-Holland,
1977.

\bibitem{Zhou_05_direction vorticity}
Y, Zhou,
\emph{A new regularity criterion for the Navier-Stokes equations in terms of the direction of vorticity,}
Monatsh. Math.
{\bf 144}
(2005),
251--257.

\bibitem{Zhou_05_one component}
Y. Zhou,
\emph{A new regularity criterion for weak solutions to the Navier-Stokes equations,}
J. Math. Pures Appl.
{\bf 84}
(2005),
1496--1514.

\bibitem{Zhou_06_grad pressure}
Y. Zhou,
\emph{On a regularity criterion in terms of the gradient of pressure for the Navier-Stokes equations in $\R^n$,}
Z. Angew. Math. Phys.
{\bf 57}
(2006),
384--392.

\bibitem{Zhou_Pokorny_09_grad one component}
Y. Zhou, M. Pokorn\'y,
\emph{On a regularity criterion for the Navier-Stokes equations involving gradient of one velocity component,}
J. Math. Phys.
{\bf 50}
(2009),
123514.

\bibitem{Zhou_Pokorny_09_one component}
Y. Zhou, M. Pokorn\'y,
\emph{On the regularity to the solutions of the Navier-Stokes equations via one velocity component,}
Nonlinearity
{\bf 23}
(2010),
1097--1107.

\bibitem{Zhou_04_pressure}
Y. Zhou,
\emph{Regularity criteria in terms of pressure for the $3$D Navier-Stokes equations in a generic domain,}
Math. Ann.
{\bf 328}
(2004),
173--192.

\bibitem{Zhou_09_weighted reg}
Y. Zhou,
\emph{Weighted regularity criteria for the three-dimensional Navier-Stokes equations,}
Proc. Roy. Soc. Edinburgh, Sect. A: Math.
{\bf 139}
(2009),
661--671.

\bibitem{Zhang_Xicheng_08}
X.C. Zhang,
\emph{A regularity criterion for the solutions of $3$D Navier-Stokes equations,}
J. Math. Anal. Appl.
{\bf 346}
(2008),
336--339.

\bibitem{Zhang_Chen_05}
Z.F. Zhang, Q.L. Chen,
\emph{Regularity criterion via two components of vorticity on weak solutions to the Navier-Stokes equations in $\R^3$,}
J. Differential Equations
{\bf 216}
(2005),
470--481.
\end{thebibliography}
\end{document}